\documentclass[12pt]{article}
\usepackage{mathrsfs}
\usepackage{stmaryrd}
\usepackage{amsfonts,amsmath,amssymb,amscd}
\usepackage{shadow}
\usepackage{hyperref}

\allowdisplaybreaks

\newtheorem{thm}{Theorem}[section]
\newtheorem{lem}[thm]{Lemma}

\newtheorem{cor}[thm]{Corollary}

\newtheorem{conj}{Conjecture}[section]

\def\qed{\hfill \rule{4pt}{7pt}}

\def\pf{\noindent {\it{Proof.} \hskip 2pt}}

\numberwithin{equation}{section}

\begin{document}
\begin{center}
{\large\bf On a Conjecture of Chen-Guo-Wang}
\date{}
\end{center}
\begin{center}
Bo Ning$^{1}$~~~~ Ken Y. Zheng$^{2}$\\[8pt]
$^{1}$Center for Applied Mathematics\\
Tianjin  University, Tianjin 300072, P. R. China\\[6pt]
$^{2}$Center for Combinatorics, LPMC-TJKLC\\
Nankai University, Tianjin 300071, P. R. China\\[6pt]
Email:$^{1}${\tt ningbo-maths@163.com}, $^{2}${\tt kenzheng@aliyun.com}
\end{center}

\begin{abstract}
Towards confirming Sun's conjecture on the strict log-concavity of combinatorial sequence involving the n$th$ Bernoulli number, Chen, Guo and Wang proposed a conjecture about the log-concavity of the function $\theta(x)=\sqrt[x]{2\zeta(x)\Gamma(x+1)}$ for $x\in (6,\infty)$, where $\zeta(x)$ is the Riemann zeta function and $\Gamma(x)$ is the Gamma function. In this paper, we first prove this conjecture along the spirit of Zhu's previous work. Second, we extend Chen et al.'s conjecture in the sense of almost infinite log-monotonicity of combinatorial sequences, which was also introduced by Chen et al. Furthermore, by using an analogue criterion to the one of Chen, Guo and Wang, we deduce the almost infinite log-monotonicity of the sequences $\frac{1}{\sqrt[n]{|B_{2n}|}}$, $T_n$ and $\frac{1}{\sqrt[n]{T_n}}$, where $B_{2n}$ and $T_{n}$ are the $2n$th Bernoulli number and the $n$th tangent number, respectively. These results can be seen as extensions of some solved conjectures of Sun.
\end{abstract}


\noindent{\bf AMS Classfication:} 05A20, 11B68

\section{Introduction}
The motivation of this article is to give a complete proof for a conjecture on the log-concavity of a special function involving the Riemann zeta function and the Gamma function, which was proposed by Chen, Guo and Wang \cite{ChenGuoWang_PEMS}, and was almost proved by Zhu \cite{Zhu2013} before.

In number theory, Firoozbakht's conjecture \cite[pp.185]{Ribenboim} states that $\sqrt[n]{p_n}>\sqrt[n+1]{p_{n+1}}$ for all $n\in \mathcal{N}$, where $p_n$ denotes the $n$th prime.
Motivated by Firoozbakht's conjecture, Sun \cite{Sun} posed a series of conjectures on the strict log-concavity of sequences in combinatorics and number theory. Recall that for a sequence $\{a_n\}_{n\geq 1}$, the sequence is called
\emph{log-concave} (resp. \emph{log-convex}) if $a^2_{n+1}\geq a_{n}a_{n+2}$ (resp. $a^2_{n+1}\leq a_{n}a_{n+2}$) for all natural number $n=1,2,\ldots$; and is called \emph{strictly log-concave} (resp. \emph{strictly log-convex}) if the corresponding inequality is strict.
Sun's conjectures stimulated a large number of further works on this topic, see \cite{ChenGuoWang,ChenGuoWang_PEMS,HouSunWen,LucaStanica,WangZhu,XiaYao,Zhao,Zhu2013}.

The Bernoulli numbers $B_0,B_1,B_2,\ldots$ are given by
$$
B_0=1,~\mbox{and}~\sum_{k=0}^{n}\binom{n+1}{k}B_k=0~~\mbox{for}~n\in \mathcal{Z}^+.
$$
It is well known that $B_k=0$ for odd $k$ and
$$
\frac{x}{e^x-1}=\sum_{n=0}^{\infty}B_n\frac{x^n}{n!},~\mbox{where}~|x|<2\pi,
$$
see \cite[pp.12]{AndrewsAskeyRoy}.

The tangent numbers
\begin{eqnarray*}
\{T(n)\}_{n\geq 0}=\{1,2,16,272,7936,\ldots\},
\end{eqnarray*}
are defined by
\begin{eqnarray*}
\tan x=\sum_{n\geq 1}T(n)\frac{x^{2n-1}}{(2n-1)!}
\end{eqnarray*}
and closely related to the Bernoulli numbers, e.g.,
\begin{eqnarray*}
T(n)=|B_{2n}|\frac{(4^n-1)4^n}{2n},
\end{eqnarray*}
and
\begin{eqnarray*}
\sqrt[n]{T(n)}=4\sqrt[n]{|B_{2n}|}\sqrt[n]{4^n-1}\sqrt[n]{\frac{1}{2n}}.
\end{eqnarray*}

In particular, Sun \cite{Sun} proposed the following conjecture on some properties of Bernoulli numbers.

\begin{conj}[Sun {\cite[Conjecture 2.15]{Sun}}]\label{ConjSun}
(i) The sequence $\{\sqrt[n]{|B_{2n}|}\}_{n\geq 1}$ is strictly increasing.
(ii) The sequence $\left\{\frac{\sqrt[n+1]{|B_{2n+2}|}}{\sqrt[n]{|B_{2n}|}}\right\}_{n\geq 2}$ is strictly decreasing.
\end{conj}

Notice that Conjecture \ref{ConjSun} (ii) states that $\{\sqrt[n]{|B_{2n}|}\}_{n\geq 2}$ is strictly log-concave, and it
is stronger than Conjecture \ref{ConjSun} (i). Luca and St\u{a}nic\u{a} \cite{LucaStanica} proved that $\{\sqrt[n]{|B_{2n}|}\}_{n\geq 2}$ is log-concave. Chen, Guo and Wang \cite{ChenGuoWang_PEMS} confirmed Conjecture \ref{ConjSun} (i) by an analytical method. Towards confirming Conjecture \ref{ConjSun} (ii), Chen, Guo and Wang \cite{ChenGuoWang_PEMS} introduced the function $\theta(x)=\sqrt[x]{2\zeta(x)\Gamma(x+1)}$, where $\zeta(x)=\sum_{n=1}^{\infty}\frac{1}{n^x}$ is the Riemann zeta function and $\Gamma(x)$ is the Gamma function, and conjectured the function is log-concave on $(6,\infty)$. As revealed by Chen et al. \cite{ChenGuoWang_PEMS}, from the well known formula $\zeta(2n)=\frac{2^{2n-1}\pi^{2n}}{(2n)!}|B_{2n}|$, one can see that $\sqrt[n]{|B_{2n}|}=\frac{1}{4\pi^{2}}\theta^2(2n)$. In view of this relationship, Chen et al. \cite{ChenGuoWang_PEMS} pointed out that the confirmation of this conjecture results in a positive answer to Conjecture \ref{ConjSun} (ii) of Sun.

\begin{conj}[Chen, Guo and Wang \cite{ChenGuoWang_PEMS}]\label{ConjCGW}
The function $\theta(x)=\sqrt[x]{2\zeta(x)\Gamma(x+1)}$ is log-concave when $x\in (6,\infty)$.
\end{conj}

A function $f$ is said to be \emph{strictly completely monotonic} on an interval $I$ if $f$ has derivatives of all orders on $I$
and $(-1)^nf^{(n)}(x)>0$ for $x\in I$ and all integers $n\geq 0$. Zhu \cite{Zhu2013} creatively applied Alzer's result \cite{Alzer}, which is concerning strict completely monotonicity of inequality involving Gamma function, to obtain the first improvement on this conjecture. He proved Conjecture \ref{ConjCGW} holds on $(7.1,\infty)$. In fact, Zhu used a general inequality for the Gamma function\cite[Lemma 3.1]{Zhu2013} which enable him to get a general result(\cite[Theorem 4.6]{Zhu2013}).

By using more precise estimates (including Lemma \ref{lem23} and the famous inequality involving the Gamma function), we confirm Chen et al.'s conjecture. Further, to know more information about log-behavior of the sequence $\{\sqrt[n]{|B_{2n}|}\}_{n\geq 1}$, we extend Chen et al.'s conjecture. Recall that in \cite{ChenGuoWang_PEMS}, Chen, Guo and Wang introduced the concept of infinite log-monotonicity of combinatorial sequences as follows:
Define an operator $R$ on a sequence $\{z_n\}_{n\geq 0}$ by  $R\{z_n\}_{n\geq 0}=\{x_n\}_{n\geq 0}$, where $x_n=z_{n+1}/z_n$. The sequence $\{z_n\}_{n\geq 0}$ is said to be \emph{infinitely log-monotonic} if the sequence $R^{r}\{z_n\}_{n\geq 0}=\{x_n\}_{n\geq 0}$ is log-concave for all odd natural number $r$ and is log-convex for all even natural number $r$. A sequence $\{a_n\}_{n\geq 0}$ is called \emph{almost infinitely log-monotonic} if for $k\geq 0$, $\{a_n\}_{n\geq 0}$ is log-monotonic of order $k$ except for certain terms at the beginning.  Chen et al. \cite{ChenGuoWang} also revealed the relationship between logarithmically completely monotonic function and infinite log-monotonicity of combinatorial sequences.

In this article, we try to extend Chen et al.'s conjecture. To prove our main result that  $(-1)^{k}[\log\theta(x)]^{(k)}<0$ where $k\geq 2$  and $x$ is large enough, we give a lower and upper bound for the $j$th derivative of $\log\Gamma(x)$, then bound the $k$th derivative of $\log\Gamma(x)\cdot\frac{1}{x}$ and the $j$th derivative of $\log\zeta(x)$, and finally consider the asymptotic property of  $\frac{x^{k+1}}{k!}(\frac{\log\zeta(x)}{x})^{(k)}$ when $x$ tends to infinity. The main theorem has some applications. By applying an analogue to a criterion of Chen, Guo and Wang (see Theorem 2.1 in \cite{ChenGuoWang}), we deduce the almost infinite log-monotonicity of the sequences $\frac{1}{\sqrt[n]{|B_{2n}|}}$, $T_n$ and $\frac{1}{\sqrt[n]{T_n}}$, respectively.

 This article is organized as follows.  In Section 2, we confirm Chen et al.'s conjecture completely in the spirit of Zhu's previous work.
In Section 3, we extend Chen et al's conjecture and obtain the almost infinite log-monotonicity of the sequences $\frac{1}{\sqrt[n]{|B_{2n}|}}$, $T_n$ and $\frac{1}{\sqrt[n]{T_n}}$.

\section{Proof of Chen-Guo-Wang's conjecture}
This section is devoted to a complete proof of Chen, Guo and Wang's conjecture.

Recall that Zhu  has used the following theorem due to Alzer \cite{Alzer}, which plays a central role in his proof.
\begin{thm}[Alzer \cite{Alzer}]
(i) The function
\begin{align}\label{G0}
G_0(x)=-\log\Gamma(x)+(x-1/2)\log x-x+\log\sqrt{2\pi}+\frac{1}{12x}
\end{align}
is strictly completely monotonic on $(0,\infty)$.\\
(ii) The function
\begin{eqnarray}\label{F0}
F_0(x)=\log\Gamma(x)-(x-1/2)\log x+x-\log\sqrt{2\pi}
\end{eqnarray}
is strictly completely monotonic on $(0,\infty)$.
\end{thm}
As indicated in \cite{Zhu2013}, the following inequalities can be easily deduced from Alzer's result:
\begin{eqnarray}\label{logG}
(x-1/2)\log x-x+\log\sqrt{2\pi}<\log\Gamma(x)<(x-1/2)\log x-x+\log\sqrt{2\pi}+\frac{1}{12x},
\end{eqnarray}
\begin{eqnarray}\label{logG'}
\log x-\frac{1}{2x}-\frac{1}{12x^2}<(\log\Gamma(x))'<\log x-\frac{1}{2x},
\end{eqnarray}
\begin{eqnarray}\label{logG''}
\frac{1}{x}+\frac{1}{2x^2}<(\log\Gamma(x))^{''}<\frac{1}{x}+\frac{1}{2x^2}+\frac{1}{6x^3}.
\end{eqnarray}

To confirm Conjecture \ref{ConjCGW}, besides more precise inequalities, we need two lemmas. The first one slightly refines Lemma 3.2 in \cite{Zhu2013}.
The second one is well known, see Andrews, Askey and Roy \cite[pp.3]{AndrewsAskeyRoy}.

\begin{lem}\label{lemRZF}
Let $\zeta(x)=\sum_{n\geq 1}\frac{1}{n^x}$ be the Riemann zeta function. Then for $x\geq 5$,
\begin{eqnarray}
1<\zeta(x)<1+\frac{1.5}{2^x}.
\end{eqnarray}
\end{lem}

\pf
The left hand is obvious. Now we will prove the right hand. Note that
\begin{eqnarray*}
\zeta(x)&=&1+\frac{1}{2^x}+\frac{1}{3^x}\left(1+\frac{1}{(\frac{4}{3})^x}+\frac{1}{(\frac{5}{3})^x}+\frac{1}{(\frac{6}{3})^x}+\ldots\right)\\[6pt]
&<&1+\frac{1}{2^x}+\frac{1}{3^x}\cdot3\cdot\left(1+\frac{1}{2^x}+\frac{1}{3^x}+\ldots\right)\\[6pt]
&=&1+\frac{1}{2^x}+\frac{1}{3^{x-1}}\zeta(x).
\end{eqnarray*}
So we have \[\zeta(x)<1+\frac{1+\frac{3^{x-1}}{2^x}}{3^{x-1}-1},\] and we only suffice to show that \[\frac{1+\frac{3^{x-1}}{2^x}}{3^{x-1}-1}<\frac{1.5}{2^x},\] that is $2^{x+1}+3<3^{x-1}$, which holds when $x\geq5$. The proof is complete.
\qed

\begin{lem}\label{lem23}
$\Gamma(x+1)=x\Gamma(x).$
\end{lem}

We are ready to give a proof of Conjecture \ref{ConjCGW}.

\begin{thm}
The function $\theta(x)=\sqrt[x]{2\zeta(x)\Gamma(x+1)}$ is log-concave when $x\in (6,\infty)$.
\end{thm}
\pf
We prove the theorem by refining the technique of Zhu \cite{Zhu2013}.
To prove Conjecture \ref{ConjCGW}, we only suffice to show that for $x\in (6,\infty)$,
$$(\log\sqrt[x]{2\zeta(x)\Gamma(x+1)})^{''}<0,$$ that is to show that for $x\in (6,\infty)$, there holds

\begin{align*}
\left(\frac{\log2+\log\zeta(x)+\log\Gamma(x+1)}{x}\right)^{''}=\left(\frac{\log2}{x}\right)^{''}+\left(\frac{\log\zeta(x)}{x}\right)^{''}+\left(\frac{\log\Gamma(x+1)}{x}\right)^{''}<0.
\end{align*}
Since
\begin{align*}
\left(\frac{\log2}{x}\right)^{''}=\frac{2\log2}{x^3},
\end{align*}
 we have
\begin{eqnarray}\label{eqn1}
x^3\left(\frac{\log2}{x}\right)^{''}=2\log2=1.386\cdots.\label{log2}
\end{eqnarray}
By computation,
\begin{eqnarray}
\left(\frac{\log\zeta(x)}{x}\right)^{''}&=&\frac{(\log\zeta(x))^{''}}{x}+2(\log\zeta(x))^{'}\left(\frac{1}{x}\right)^{'}+\log\zeta(x)\left(\frac{1}{x}\right)^{''}\nonumber\\[6pt]
&=&\frac{{\zeta(x)}^{''}\zeta(x)-{{\zeta(x)}^{'}}^{2}}{{\zeta(x)}^{2}}\cdot\frac{1}{x}-\frac{2}{x^2}\frac{{\zeta(x)}^{'}}{\zeta(x)}+2\log\zeta(x)\cdot\frac{1}{x^3}\nonumber\\[6pt]
&\leq & {\zeta(x)}^{''}\cdot\frac{1}{x}+2\frac{|{\zeta(x)}^{'}|}{x^2}+\frac{2\log\zeta(x)}{x^3}.\label{logzeta}
\end{eqnarray}

According to Lemma \ref{lemRZF}, noted that $1+\frac{1.5}{2^x}$ is decreasing and $1+\frac{1.5}{2^5}<1.047$, we have that for $x\geq5$, $\zeta(x)\in(1,1.047)$.
It is easy to get that $\log x-\sqrt{x}$ is increasing when $x\in [1,1.047]$, then\[\log\zeta(x)-\sqrt{\zeta(x)}< \log1.047-\sqrt{1.047}<-0.977.\]
Thus
\begin{eqnarray}
\log\zeta(x)< \sqrt{\zeta(x)}-0.977<\sqrt{1+\frac{1.5}{2^x}}-0.977.\label{zeta}
\end{eqnarray}
By Lemma \ref{lemRZF}, we have
\begin{eqnarray}
{\zeta(x)}^{''}=\sum_{n\geq 2}\frac{(\log n)^2}{n^x}<\sum_{n\geq 2}\frac{1}{n^{x-1}}=\zeta(x-1)-1\leq \frac{1.5}{2^{x-1}},\label{zeta2}
\end{eqnarray}
and
\begin{eqnarray}
|{\zeta(x)}^{'}|=\sum_{n\geq 2}\frac{(\log n)}{n^x}<\sum_{n\geq 2}\frac{1}{n^{x-\frac{1}{2}}}=\zeta(x-\frac{1}{2})-1\leq \frac{1.5}{2^{x-\frac{1}{2}}}.\label{zeta1}
\end{eqnarray}
So putting inequalities (\ref{zeta2}), (\ref{zeta1}) and (\ref{zeta}) into (\ref{logzeta}), respectively, by Lemma \ref{lemRZF}, we have
\begin{eqnarray}
x^3\left(\frac{\log\zeta(x)}{x}\right)^{''}&<\frac{1.5x^2}{2^{x-1}}+\frac{3x}{2^{x-\frac{1}{2}}}+2(\sqrt{1+\frac{1.5}{2^x}}-0.977). \label{zeta2x3}
\end{eqnarray}
Define $$f_0(x)=\frac{1.5}{2^{x-1}}(x^2+\sqrt{2}x).$$ By derivation, we can get
\begin{eqnarray}
f_0^{'}(x)=-\frac{3}{2^x}(\log2\cdot x^2+(\sqrt{2}\log2-2)x-\sqrt{2})<0 \label{f01}
\end{eqnarray}
for $x\geq 3$. Together with (\ref{zeta2x3}) and (\ref{f01}), we know
\begin{eqnarray}
x^3\left(\frac{\log\zeta(x)}{x}\right)^{''}<f_0(6)+2\left(\sqrt{1+\frac{1.5}{2^6}}-0.977\right)=2.1545\cdots\label{x3}
\end{eqnarray}
when $x\in (6,\infty)$. Furthermore, we have $$\left(\frac{\log\Gamma(x+1)}{x}\right)^{''}=\left(\frac{\log x}{x}\right)^{''}+\left(\frac{\log\Gamma(x)}{x}\right)^{''},$$
where
\begin{eqnarray}
\left(\frac{\log x}{x}\right)^{''}=-\frac{3}{x^3}+\frac{2\log x}{x^3}.\label{logx/x}
\end{eqnarray}
Note that
\begin{eqnarray}
\left(\frac{\log\Gamma(x)}{x}\right)^{''}=\frac{(\log\Gamma(x))^{''}}{x}-\frac{2(\log\Gamma(x))^{'}}{x^2}+\frac{2\log\Gamma(x)}{x^3}.\label{logzetax}
\end{eqnarray}
Applying inequalities \eqref{logG}, \eqref{logG'} and \eqref{logG''} to the terms $(\log\Gamma(x))^{''}$, $(\log\Gamma(x))^{'}$ and $\log\Gamma(x)$ in (\ref{logzetax}), respectively, we obtain
\begin{eqnarray}
\left(\frac{\log\Gamma(x)}{x}\right)^{''}&<&\frac{1}{x}\left(\frac{1}{x}+\frac{1}{2x^2}+\frac{1}{6x^3}\right)-\frac{2}{x^2}\left(\log x-\frac{1}{2x}-\frac{1}{12x^2}\right)\nonumber\\[6pt]
&&+\frac{2}{x^3}\left((x-\frac{1}{2})\log x-x+\log\sqrt{2\pi}+\frac{1}{12x}\right)\nonumber\\[6pt]
&=&-\frac{1}{x^2}+\frac{3}{2x^3}-\frac{\log x}{x^3}+\frac{\log2\pi}{x^3}+\frac{1}{2x^4}.\label{Gammax/x}
\end{eqnarray}
Combining (\ref{logx/x}) and (\ref{Gammax/x}), by Lemma \ref{lemRZF}, we have
\begin{eqnarray*}
\left(\frac{\log\Gamma(x+1)}{x}\right)^{''}<-\frac{1}{x^2}-\frac{3}{2x^3}+\frac{\log x}{x^3}+\frac{\log2\pi}{x^3}+\frac{1}{2x^4}.
\end{eqnarray*}
Furthermore, we have
\begin{eqnarray*}
x^3\left(\frac{\log\Gamma(x+1)}{x}\right)^{''}<-x+\log x-\frac{3}{2}+\frac{1}{2x}+\log2\pi.
\end{eqnarray*}
Define \[f_1(x)=-x+\log x-\frac{3}{2}+\frac{1}{2x}+\log2\pi.\] Since the derivation \[f'_1(x)=-1+\frac{1}{x}-\frac{1}{2x^2}=-\frac{1}{2x^2}(2x^2-2x+1)<0,\] we obtain
\begin{eqnarray}
x^3\left(\frac{\log\Gamma(x+1)}{x}\right)^{''}<f_1(6)=-6+\log6-\frac{3}{2}+\frac{1}{12}+\log2\pi=-3.787\cdots\label{Gammax+1}.
\end{eqnarray}
By (\ref{log2}), (\ref{x3}) and (\ref{Gammax+1}), we can obtain $$\left(\frac{\log2+\log\zeta(x)+\log\Gamma(x+1)}{x}\right)^{''}<-0.2465<0.$$ The proof is complete.
\qed

\section{An extension of Chen et al.'s conjecture and applications}
 In this section, we extend Chen et al.'s conjecture. By proving that  $(-1)^{k}[\log\theta(x)]^{(k)}<0$ when $k\geq 2$  and $x$ is large enough, we prove the almost infinitely log-monotonic property of three sequences, including $\frac{1}{\sqrt[n]{|B_{2n}|}}$, $T_n$ and $\frac{1}{\sqrt[n]{T_n}}$. Our main theorem of this section is the following.

\begin{thm}\label{thm-th}
Let $\theta(x)=\sqrt[x]{2\zeta(x)\Gamma(x+1)}$. Then there holds $(-1)^{k}[\log\theta(x)]^{(k)}<0$, where $k\geq 2$ and $x$ is large enough.
\end{thm}

To prove Theorem \ref{thm-th}, we need several lemmas. The first two lemmas are quite easy, so we only give the proof of Lemma \ref{b-jthlogG}.

\begin{lem}
For $j\geq1$,
\begin{eqnarray}\label{jth-logx/x}
\left(\frac{\log x}{x}\right)^{(j)}=(-1)^{j-1}j!\left(\sum_{i=1}^{j}\frac{1}{i}\cdot\frac{1}{x^{j+1}}-\frac{\log x}{x^{j+1}}\right).
\end{eqnarray}
\end{lem}

\begin{lem}
(i) For even $j\geq 2$, $x\in (0,\infty)$,
\begin{eqnarray}\label{jth-logG-e}
0<j!\left(\frac{1}{j(j\!-\!1)x^{j-1}}\!+\!\frac{1}{2jx^j}\right)<[\log\Gamma(x)]^{(j)}<j!\left(\frac{1}{j(j-1)x^{j\!-\!1}}\!+\!\frac{1}{2jx^j}\!+\!\frac{1}{12x^{j+1}}\right).
\end{eqnarray}
(ii) For odd $j\geq 3$, $x\in (0,\infty)$,
\begin{eqnarray}\label{jth-logG-o}
-\!j!\left(\frac{1}{j(j\!-\!1)x^{j-1}}\!+\!\frac{1}{2jx^j}\!+\!\frac{1}{12x^{j+1}}\right)\!<\![\log\Gamma(x)]^{(j)}\!<\!-j!\left(\frac{1}{j(j\!-\!1)x^{j-1}}\!+\!\frac{1}{2jx^j}\right)<0.
\end{eqnarray}
\end{lem}

\begin{lem}\label{b-jthlogG}
(i) For even $k\geq 2$,
\begin{eqnarray}
\left(\log\Gamma(x)\cdot\frac{1}{x}\right)^{(k)}<\frac{k!}{x^{k+1}}\left(-\frac{x}{k}-\frac{\log x}{2}+\frac{\log\pi}{2}+\frac{1}{2}\sum_{i=1}^k\frac{1}{i}+\frac{k+1}{12x}\right).
\end{eqnarray}
(ii) For odd $k\geq 3$,
\begin{eqnarray}
\left(\log\Gamma(x)\cdot\frac{1}{x}\right)^{(k)}>-\frac{k!}{x^{k+1}}\left(-\frac{x}{k}-\frac{\log x}{2}+\frac{\log\pi}{2}+\frac{1}{2}\sum_{i=1}^k\frac{1}{i}+\frac{k+1}{12x}\right).
\end{eqnarray}
\end{lem}

\pf
(i) By Lebniz formula,
\begin{eqnarray*}
\left(\log\Gamma(x)\cdot\frac{1}{x}\right)^{(k)}=\sum_{j=0}^{k}\binom{k}{j}(\log\Gamma(x))^{(j)}\left(\frac{1}{x}\right)^{(k-j)}.
\end{eqnarray*}
Let $k$ be an even integer and $k\geq2$. The cases of $j=0$ and $j=1$ can be easily obtained from \eqref{logG} and \eqref{logG'}, respectively.
Now we consider the case when $2\leq j\leq k$. If $j$ is even, we have
\[(\frac{1}{x})^{(k-j)}=(-1)^{k-j}\frac{(k-j)!}{x^{k-j+1}}=\frac{(k-j)!}{x^{k-j+1}}>0.\]
By \eqref{jth-logx/x}, we see that
\begin{eqnarray*}
0<[\log\Gamma(x)]^{(j)}<j!\left(\frac{1}{j(j-1)x^{j-1}}+\frac{1}{2jx^j}+\frac{1}{12x^{j+1}}\right).
\end{eqnarray*}
So
\begin{eqnarray}
\notag\binom{k}{j}(\log\Gamma(x))^{(j)}(\frac{1}{x})^{(k-j)}&<&\frac{k!}{j!(k-j)!}j!\left(\frac{1}{j(j-1)x^{j-1}}+\frac{1}{2jx^j}+\frac{1}{12x^{j+1}}\right)\cdot\frac{(k-j)!}{x^{k-j+1}}\nonumber\\[6pt]
&=&\frac{k!}{x^{k+1}}\left(\frac{x}{j(j-1)}+\frac{1}{2j}+\frac{1}{12x}\right).\label{binsum}
\end{eqnarray}
If $j$ is odd, we have
\[\left(\frac{1}{x}\right)^{(k-j)}=(-1)^{k-j}\frac{(k-j)!}{x^{k-j+1}}=\frac{(k-j)!}{x^{k-j+1}}<0.\]
By \eqref{jth-logG-e}, we get
\begin{eqnarray*}
-j!\left(\frac{1}{j(j-1)x^{j-1}}+\frac{1}{2jx^j}+\frac{1}{12x^{j+1}}\right)<[\log\Gamma(x)]^{(j)}<0.
\end{eqnarray*}
So \eqref{binsum} also holds in this case. Now we obtain
\begin{eqnarray*}
&&\sum_{j=0}^{k}\binom{k}{j}(\log\Gamma(x))^{(j)}\left(\frac{1}{x}\right)^{(k-j)}\nonumber\\[6pt]
&<&\log\Gamma(x)\cdot(-1)^kk!\frac{1}{x^{k+1}}+k(\log\Gamma(x))^{'}\cdot(-1)^{k-1}\frac{(k-1)!}{x^k}\nonumber\\[6pt]
&&+\sum_{j=1}^{k}\frac{k!}{x^{k+1}}\left(\frac{x}{j(j-1)}+\frac{1}{2j}+\frac{1}{12x}\right)\\[6pt]
&<&\frac{k!}{x^{k+1}}\left((x-1/2)\log x-x+\log\sqrt{2\pi}+\frac{1}{12x}\right)+k!\frac{1}{x^{k+1}}\left(-x\log x+\frac{1}{2}+\frac{1}{12x}\right)\nonumber\\[6pt]
&&+k!\left[\left(1-\frac{1}{k}\right)x+\frac{1}{2}\sum_{j=2}^k\frac{1}{j}+\frac{k-1}{12x}\right]\nonumber\\[6pt]
&=&\frac{k!}{x^{k+1}}\left(-\frac{x}{k}-\frac{\log x}{2}+\frac{\log2\pi}{2}+\frac{1}{2}\sum_{j=1}^{k}\frac{1}{j}+\frac{k+1}{12x}\right).
\end{eqnarray*}
(ii) The case of odd $k$ when $k\geq3$ is similar, we omit the proof.
\qed

\begin{lem}\label{b-jthlogz}
Let $k\geq 2$ be an integer and $0\leq j\leq k$. If $k$ is even, then
\begin{eqnarray}
(-1)^{k-j}\log^{(j)}\zeta(x)\leq \frac{1.5}{2^{x-\frac{j}{2}}};
\end{eqnarray}
If $k$ is odd, then
\begin{eqnarray}
0>(-1)^{k-j}\log^{(j)}\zeta(x)\geq -\frac{1.5}{2^{x-\frac{j}{2}}}.
\end{eqnarray}
\end{lem}

\pf
First assume that $k$ is even. Recall that
\begin{eqnarray*}
-\frac{\zeta^{'}(x)}{\zeta(x)}=\sum_{n=1}^{\infty}\frac{\Lambda(n)}{n^x},
\end{eqnarray*}
where $\Lambda(n)$ is the von Mangoldt function, see \cite[pp.122]{AndrewsAskeyRoy}. So if $j\geq 1$ then
\begin{eqnarray}
\notag (-1)^{k-j}\log^{(j)}\zeta(x)&=&(-1)^{k}\sum_{n=2}^{\infty}\frac{\Lambda(n)(\log n)^{j-1}}{n^x}\\[6pt]
&\leq&\sum_{n=2}^{\infty}\frac{(\log n)^{j}}{n^x}.\label{jth-logz}
\end{eqnarray}
Since $\log x\leq \sqrt{x}-1$ when $x\in [1,\infty)$, \eqref{jth-logz} and Lemma \ref{lemRZF} imply that
\begin{eqnarray*}
(-1)^{k-j}\log^{(j)}\zeta(x)&\leq& \sum_{n=2}^{\infty}\frac{1}{n^{x-\frac{1}{2}}}\leq \frac{1.5}{2^{x-\frac{j}{2}}}.
\end{eqnarray*}
If $j=0$, then $(-1)^k\log\zeta(x)\leq \zeta(x)\leq \frac{1.5}{2^{x}}$ by Lemma \ref{lemRZF}.

The case when $k$ is odd can be proved similarly. We omit the details.
\qed

\begin{lem}
For any $k\geq 2$,
\begin{eqnarray}
\frac{x^{k+1}}{k!}\left(\frac{\log\zeta(x)}{x}\right)^{(k)}=o(1)
\end{eqnarray}
\end{lem}

\pf
If $k$ is even, then we have
\begin{align*}
x^{k+1}\left(\frac{\log\zeta(x)}{x}\right)^{(k)}>0
\end{align*}
 and
\begin{eqnarray}
\left(\frac{\log\zeta(x)}{x}\right)^{(k)}&=&\sum_{j=0}^{k}\binom{n}{k}[\log\zeta(x)]^{(j)}\left(\frac{1}{x}\right)^{(k-j)}\nonumber\\[6pt]
&=&\sum_{j=0}^k\frac{k!}{j!(k-j)!}(-1)^{k-j}[\log\zeta(x)]^{(j)}(k-j)!\frac{1}{x^{k-j+1}}\nonumber\\[6pt]
\label{in1}&\leq& k!\sum_{j=0}^k\frac{1}{j!}\cdot\frac{1.5}{2^{x-\frac{j}{2}}}\cdot\frac{1}{x^{k-j+1}}\\[6pt]
\label{in2}&\leq& \frac{k!}{x^{k+1}}\cdot 1.5e\cdot\frac{\sum_{j=0}^k(\sqrt{2}x)^j}{2^x},
\end{eqnarray}
where \eqref{in1} holds by Lemma \ref{b-jthlogz}, and \eqref{in2} holds by considering the equality $e=\sum_{j=0}^{\infty}\frac{1}{j!}$. Similarly, if $k\geq 3$ and $k$ is odd, then
we have \[-\frac{1}{x^{k+1}}\cdot 1.5e\cdot\frac{(\sqrt{2}x)^k}{2^x}<\frac{x^{k+1}}{k!}\left(\frac{\log\zeta(x)}{x}\right)^{(k)}<0.\]
Note that \[\lim_{x\rightarrow\infty}k!\cdot 1.5e\cdot\frac{(\sqrt{2}x)^k}{2^x}=0.\]
So
\begin{align*}
\ \ \ \ \ \ \ \ \ \ \ \ \ \ \ \ \ \ \ \ \ \ \ \ \ \ \ \ \ \ \ \ \ \ \ \ \ \ \ x^{k+1}\cdot \left(\frac{\log\zeta(x)}{x}\right)^{(k)}=o(1).\ \ \ \ \ \ \ \ \ \ \ \ \ \ \ \ \ \ \ \ \ \ \ \ \ \ \ \ \ \ \ \ \ \ \ \ \ \ \ \qed
\end{align*}

Now we give a proof of Theorem \ref{thm-th}.

\noindent{}
{\bf Proof of Theorem \ref{thm-th}}. If $k$ is even, then
\begin{align*}
 &(-1)^{k}x^{k+1}[\log\theta(x)]^{(k)}\\[6pt]
 &=x^{k+1}\left(\left(\frac{\log2}{x}\right)^{(k)}+\left(\frac{\log\zeta(x)}{x}\right)^{(k)}+\left(\frac{\log\Gamma(x+1)}{x}\right)^{(k)}\right)\\[6pt]
 &\leq k!\left(\log2+o(1)+\log x-\sum_{j=1}^{k}\frac{1}{j}-\frac{x}{k}+\frac{1}{2}\sum_{j=1}^{k}\frac{1}{j}+\log\sqrt{2\pi}-\frac{\log x}{2}+\frac{k+1}{12x}\right)\\[6pt]
&= k!\left(\frac{\log x}{2}-\frac{1}{2}\sum_{j=1}^{k}\frac{1}{j}-\frac{x}{k}+\frac{3\log2}{2}+\frac{\log\pi}{2}+\frac{k+1}{12x}+o(1)\right).
\end{align*}
Define
\begin{align*}
f(k,x)=\frac{\log x}{2}-\frac{1}{2}\sum_{j=1}^{k}\frac{1}{j}-\frac{x}{k}+\frac{3\log2}{2}+\frac{\log\pi}{2}+\frac{k+1}{12x}.
\end{align*}
For a given $k$, $$\frac{df}{dx}=\frac{1}{2x}-\frac{1}{k}-\frac{k+1}{12x^2}=-\frac{12x^2-6kx+k^2+1}{12kx^2}<0.$$
So we have that
\begin{align*}
\notag f(k,x)\leq f(k,3k)&=\frac{\log3k}{2}-\frac{1}{2}\sum_{j=1}^{k}\frac{1}{j}-3+\frac{3\log2}{2}+\frac{\log\pi}{2}+\frac{k+1}{36k}\\[6pt]
\notag &=\frac{\log3}{2}-\frac{1}{2k}-3+\frac{3\log2}{2}+\frac{\log\pi}{2}+\frac{1}{36}+\frac{1}{36k}\\[6pt]
&=-\frac{17}{36k}-0.8108..
\end{align*}
Hence $(-1)^{k}x^{k+1}[\log\theta(x)]^{(k)}<0$ when $x$ is large enough.

 If $k$ is odd, then
\begin{align*}
 &x^{k+1}[\log\theta(x)]^{(k)}\\[6pt]
 &=x^{k+1}\left(\left(\frac{\log2}{x}\right)^{(k)}+\left(\frac{\log\zeta(x)}{x}\right)^{(k)}+\left(\frac{\log\Gamma(x+1)}{x}\right)^{(k)}\right)\\[6pt]
&\geq k!\left(-\log2+o(1)-\log x+\sum_{j=1}^{k}\frac{1}{j}+\frac{x}{k}-\frac{1}{2}\sum_{j=1}^{k}\frac{1}{j}-\log\sqrt{2\pi}+\frac{\log x}{2}-\frac{k+1}{12x}\right)\\[6pt]
&= k!\left(-\frac{\log x}{2}+\frac{1}{2}\sum_{j=1}^{k}\frac{1}{j}+\frac{x}{k}-\frac{3\log2}{2}-\frac{\log\pi}{2}-\frac{k+1}{12x}+o(1)\right)\\[6pt]
&=k!(-f(k,x)+o(1)).
\end{align*}
As revealed in the above, we have $-f(k,x)+o(1)>0$ when $k$ is odd. So we can see that when $k$ is odd and $x$ is large enough, $(-1)^{k}x^{k+1}[\log\theta(x)]^{(k)}<0$.

Hence $(-1)^{k}x^{k+1}[\log\theta(x)]^{(k)}<0$ when $x$ is large enough.\qed

It is obviously that Theorem \ref{thm-th} has the following immediate corollary.
\begin{cor}
The function $\theta^{-1}(x)=\frac{1}{\sqrt[x]{2\zeta(x)\Gamma(x+1)}}$ is almost completely log-monotonic.
\end{cor}


In \cite{ChenGuoWang}, Chen et al. proved the following useful theorem.
\begin{thm}
Assume that $f(x)$ is a function such that $[\log f(x)]''$ is completely monotonic for $x\geq 1$.
Let $a_n=f(n)$ for $n\geq 1$. Then the sequence $\{a_n\}_{n\geq 1}$ is infinitely log-monotonic.
\end{thm}

Moreover, one could get the following result on almost infinitely log-monotonic sequences. To ensure the integrity of this paper, we include all the details here, which is similar to the proof of Theorem 2.1 in \cite{ChenGuoWang}.
\begin{thm}\label{sim-th}
Assume that $f(x)$ is a function such that $[\log f(x)]''$ is almost completely monotonic.
Let $a_n=f(n)$ for $n\geq 1$. Then the sequence $\{a_n\}_{n\geq 1}$ is almost infinitely log-monotonic.
\end{thm}
\pf We give the proof step by step as in \cite[Theorem 2.1]{ChenGuoWang}. Let $b_{n,0}=a_n$ and $b_{n,i+1}=b_{n+1,i}/b_{n,i}$. Set $f_0(x)=f(x)$, and define the functions $f_1(x)$, $f_2(x)\cdots$ by the relation
\begin{align}
f_{i+1}(x)=\frac{f_i(x+1)}{f_i(x)}.
\end{align}
Since $b_{n,i}=f_i(n)$ for any $i\geq0$ and $n\geq1$, we shall show that for $j\geq0$, $k\geq2$ and enough large $x$,
\begin{align}\label{2.2}
(-1)^k[\log f_{2j}(x)]^{(k)}\geq0,
\end{align}
and
\begin{align}\label{2.3}
(-1)^k[\log f_{2j+1}(x)]^{(k)}\leq0.
\end{align}
We prove it by induction on $j$. Since $[\log f(x)]^{''}$ is almost completely monotonic, we see that for $k\geq2$ and enough large $x$,
\begin{align}\label{2.4}
(-1)^k[\log f(x)]^{(k)}\geq0,
\end{align}
that is \eqref{2.2} holds for $j=0$. The inequality \eqref{2.4} reveals that for $k\geq1$ and enough large $x$,
\begin{align}\label{2.5}
(-1)^k[\log f(x)]^{(k+1)}\leq0.
\end{align}
Because of $f_i(x)=f(x+1)/f(x)$, by \eqref{2.5} we see that for $k\geq2$ and enough large $x$,
\[(-1)^k[\log f_1(x)]^{(k)}=(-1)^k[\log f(x+1)]^{(k)}-(-1)^k[\log f(x)]^{(k)}\leq0.\]
Thus \eqref{2.3} is true for $j=0$.

We now assume that \eqref{2.2} and \eqref{2.3} for $j\geq n-1$. We are faced to show that \eqref{2.2} and \eqref{2.3} for $j=n$. From the induction hypothesis \eqref{2.3} we see that for $k\geq2$ and enough large $x$,
\[(-1)^k[\log f_{2n}(x)]^{(k)}=(-1)^k[\log f_{2n-1}(x+1)]^{(k)}-(-1)^k[\log f_{2n-1}(x)]^{(k)}\geq0.\]
So \eqref{2.2} holds for $j=n$. Similarly, it can be shown that for $k\geq2$ and enough large $x$,
\[(-1)^k[\log f_{2n+1}(x)]^{(k)}\geq0,\]
that is, \eqref{2.3} holds for $j=n$.

Combining \eqref{2.2} and \eqref{2.3}, we conclude that for any $i\geq0$, the sequence $\{f_{2i}(n)\}_{n\geq1}$ is log-convex and the sequence $\{f_{2i+1}(n)\}_{n\geq1}$ is log-concave for enough large $n$, which completes the proof.
\qed

As a byproduct, we can obtain the almost infinite log-monotonicity of the sequence $\left\{\frac{1}{\sqrt[n]{|B_{2n}|}}\right\}_{n\geq 1}$.
\begin{cor}
The sequence $\left\{\frac{1}{\sqrt[n]{|B_{2n}|}}\right\}_{n\geq 1}$ is almost infinitely log-monotonic.
\end{cor}
\pf Recall that
\[\frac{1}{\sqrt[n]{B_{2n}}}=4\pi^2\theta^{-2}(2n).\]
Set
\begin{align*}
y(x)=4\pi^2\theta^{-2}(2x).
\end{align*}
So we have $y(n)=\frac{1}{\sqrt[n]{B_{2n}}}$. Since for $k\geq2$ and $x\geq1$,
\[(\log y(x))^{(k)}=-2(\log\theta(2x))^{(k)},\]
we see that $(-1)^{k}(\log y(x))^{(k)}>0$ for $k\geq2$ and $x$ is large enough. Then $[\log y(x)]^{''}$ is almost completely monotonic.
Thus the sequence $\left\{\frac{1}{\sqrt[n]{|B_{2n}|}}\right\}_{n\geq 1}$ is almost infinitely log-monotonic by Theorem \ref{sim-th}.
\qed

We also study the almost infinitely log-monotonic property of sequences involving tangent numbers.

\begin{thm}
The sequence $T_n$ is almost infinitely log-monotonic.
\end{thm}
\pf
As in \cite{ChenGuoWang}, set $$z(x)=\frac{2\zeta(2x)\Gamma(2x+1)}{(2\pi)^{2x}},$$ where $\zeta(x)$ is the zeta function and
$\Gamma(x)$ is the Gamma function. Let $$t(x)=z(x)\frac{(4^x-1)4^x}{2x}.$$ Then note that
$t(n)=T_n$. For $k\geq 2$,
\begin{eqnarray*}
\big(\log t(x))^{(k)}=\big(\log z(x))^{(k)}+\big(\log (4^x-1))^{(k)}+\frac{(-1)^k(k-1)!}{x^k}.
\end{eqnarray*}
In \cite{ChenGuoWang}, Chen, Guo and Wang have proved that $(-1)^k\big(\log z(x))^{(k)}>0$.

In the following, we will give a useful estimate about $(\log (4^x-1))^{(k)}$, which can be proved
by induction and simple computation. We omit the details.

{\bf Claim.} For $k\geq 2$, $x\neq 0$,
\begin{eqnarray*}
|\big(\log (4^x-1))^{(k)}|<\sum_{i=1}^{k}\frac{(\log4)^k(k-1)!}{(4^x-1)^i}.
\end{eqnarray*}

Note that
\begin{align*}
&(-1)^k\cdot\left((\log (4^x-1))^{(k)}+\frac{(-1)^k(k-1)!}{x^k}\right)\\[6pt]
&\ \ >\frac{(k-1)!}{x^k}-\sum_{i=1}^{k}\frac{(\log4)^k(k-1)!}{(4^x-1)^i}\\[6pt]
&\ \ >(k-1)!\left(\frac{1}{x^k}-\frac{(\log 4\cdot x)^k}{4^x-2}\right)\\[6pt]
&\ \ >0
\end{align*}
when $x$ is largely enough, where the last inequality holds since the simple fact that
\begin{eqnarray*}
\lim_{x\rightarrow \infty}\frac{(\log 4\cdot x)^k}{4^x-2}=0.
\end{eqnarray*}
Thus, $(-1)^k\big(\log t(x))^{(k)}>0$ when $x$ is large enough. By Theorem \ref{sim-th}, the proof is complete.
\qed

\begin{thm}
The sequence $\{\frac{1}{\sqrt[n]{T_n}}\}_{n\geq 1}$ is almost infinitely log-monotonic.
\end{thm}
\pf
 Recall that  $$\theta(x)=\sqrt[x]{2\zeta(x)\Gamma(x+1)},$$ where $\zeta(x)=\sum_{n=1}^{\infty}\frac{1}{n^x}$ is the Riemann zeta function and $\Gamma(x)$ is the Gamma function.
Now we obtain $$\sqrt[n]{|B_{2n}|}=\frac{1}{4\pi^{2}}\theta^2(2n),$$ and
\begin{eqnarray*}
\sqrt[n]{T(n)}=\frac{1}{\pi^{2}}\theta^2(2n)\sqrt[n]{4^n-1}\sqrt[n]{\frac{1}{2n}}.
\end{eqnarray*}
Thus for $k\geq 2$,
\begin{eqnarray*}
\left(\log \frac{1}{\sqrt[x]{t(x)}}\right)^{(k)}=-(\log \theta(2x))^{(k)}-\left(\frac{\log (4^x-1)}{x}\right)^{(k)}+\left(\frac{\log2}{x}\right)^{(k)}+\left(\frac{\log x}{x}\right)^{(k)}.
\end{eqnarray*}
By Theorem \ref{thm-th}, for any integer $k$ and large enough $x$,
\begin{eqnarray*}
(-1)^k(\log\theta^{-1}(2x))^{(k)}>0.
\end{eqnarray*}
 Furthermore, by induction,
one can obtain $$(-1)^{k+1}\big(\log (4^x-1))^{(k)}>0.$$ Thus, we get
\begin{align*}
&(-1)^{k+1}\left(\frac{\log (4^x-1)}{x}\right)^{(k)}\\[6pt]
&\ \ =\sum_{i=0}^{k}\binom{k}{i}(-1)^{i+1}(\log (4^x-1))^{(i)}(-1)^{2(k-i)}(k-i)!\frac{1}{x^{k-i+1}}\\[6pt]
&\ \ =\sum_{i=0}^{k}\frac{k!}{i!}[(-1)^{i+1}(\log (4^x-1))^{(i)}]\frac{1}{x^{k-i+1}}\\[6pt]
&\ \ >0.
\end{align*}
Furthermore,
\begin{eqnarray*}
(-1)^k\left(\frac{\log2}{x}\right)^{(k)}=\frac{\log2\cdot k!}{x^{k+1}};
\end{eqnarray*}
and by Lemma \ref{jth-logx/x},
\begin{eqnarray*}
(-1)^k\left(\frac{\log x}{x}\right)^{(k)}=\frac{k!}{x^{k+1}}\left(\log x-\sum_{i=1}^{k}\frac{1}{i}\right)\geq 0.
\end{eqnarray*}
The above inequality is true because we could choose $x$ which is larger than $\log k+\gamma$, where $\gamma$ is Euler-Mascheroni constant.

Thus $(-1)^k\left(\log \frac{1}{\sqrt[x]{t(x)}}\right)^{(k)}>0$ when $x$ is large enough, which completes the proof.
\qed

\noindent{\bf Remark. }Sun\cite{Sun} conjectured that $\sqrt[n]{T(n)}$ is strictly log-concave (see Sun \cite[Conjecture~3.5]{Sun}), and this conjecture was confirmed by Luca and St\u{a}nic\u{a} \cite{LucaStanica}, and Zhu \cite{Zhu2013}, independently, by different methods. Our Theorem \ref{thm-th} can be seen as
an extension of Sun's conjecture in some sense.

\noindent{\bf Added Note. } This paper is an extended version of arXiv:1508.01793, which includes a complete proof of Chen et al.'s conjecture (submitted on 3 Aug 2015).

\vspace{0.5cm}
 \noindent{\bf Acknowledgments.} This work was supported by  the 973
Project, the PCSIRT Project of the Ministry of Education,  and the National Science
Foundation of China.

\end{document}